\documentclass[12pt]{amsart}
\usepackage{a4wide}
\usepackage{eucal}
\usepackage{latexsym}
\usepackage{amssymb}
\usepackage[pdfauthor   = {Mohamed Barakat},
	    pdftitle    = {The Existence of Cartan Connections and Geometrizable Principal Bundles},
	    pdfsubject  = {53B05; 53C05; 53C10},
	    pdfkeywords = {Cartan connections; soldering form; geometrizable principal bundles},
	    bookmarks=true,
	    bookmarksopen=true,
	    colorlinks=true,
	    pagebackref=true,
	    hyperindex=true,
	    linkcolor=blue,
	    pagecolor=blue,
	    citecolor=blue,
	    urlcolor=blue,
	    ps2pdf=true,
	    ]{hyperref}

\def\bdoc{
\begin{document}}
\def\edoc{\end{document}}
\def\btp{\begin{titlepage}}
\def\etp{\end{titlepage}}
\def\bab{\begin{abstract}}
\def\eab{\end{abstract}}
\def\sec{\section}
\def\subs{\subsection}
\def\bea{\begin{array}{r@{,\ \ \ \ \ \ \ }l}}
\def\ena{\end{array}}
\def\barr{\begin{array}}
\def\earr{\end{array}}
\def\ben{\begin{enumerate}}
\def\een{\end{enumerate}}
\def\beq{\begin{equation}}
\def\eeq{\end{equation}}
\def\bqa{\begin{eqnarray*}}
\def\eqa{\end{eqnarray*}}
\def\bqan{\begin{eqnarray}}
\def\eqan{\end{eqnarray}}
\def\bit{\begin{itemize}}
\def\eit{\end{itemize}}
\def\bfi{\begin{figure}[htbp]}
\def\efi{\end{figure}}
\def\bce{\begin{center}}
\def\ece{\end{center}}
\def\bbb{\begin{thebibliography}}
\def\ebb{\end{thebibliography}}
\def\btb{\begin{tabular}}
\def\etb{\end{tabular}}
\def\btbl{\begin{table}[htbp]}
\def\etbl{\end{table}}
\def\idx{\index}

%% skips and indents:

\def\bsn{\bigskip \noindent}
\def\bs{\bigskip}
\def\br{\\[0.05cm]}
\def\ind{\indent}
\def\no{\noindent}
\def\ssn{\smallskip \noindent}
\def\ssk{\smallskip}
\def\bl{$\bullet \ $}
\def\1{{\bf 1.}}
\def\2{{\bf 2.}}

%% math abbreviations:

\def\im{\Rightarrow}
\def\ra{\rightarrow}
\def\map{\mapsto}
\def\ch{\choose}
\def\eq{\equiv}
\def\dv{\uplus}
\def\N{{\mathbb{N}}}
\def\Z{{\mathbb{Z}}}
\def\K{{\mathbb{K}}}
\def\F{{\mathbb{F}}}
\def\Q{{\mathbb{Q}}}
\def\R{{\mathbb{R}}}
\def\C{{\mathbb{C}}}
\def\H{{\mathbb{H}}}
\def\gdw{\Longleftrightarrow}
\def\sub{\subseteq}
\def\all{\forall}
\def\ex{\exists}
\def\Prf{{\sc Proof}}
\def\prf{{\sc Proof. \quad}}
\def\df{\em}
\def\inp{\ \lrcorner\ }
\def\wdg{\wedge}
\def\surj{\twoheadrightarrow}
\def\sr{\stackrel}
\def\trleq{\trianglelefteq}
\def\<{\langle}
\def\>{\rangle}

\def\nthm{\newtheorem{axiom}{Theorem}[section]}
\def\bthm{\begin{axiom}}
\def\ethm{\end{axiom}}
\def\ncn{\newtheorem{conj}[axiom]{Conjecture}}
\def\bcn{\begin{conj}}
\def\ecn{\end{conj}}
\def\nla{\newtheorem{lemm}[axiom]{Lemma}}
\def\bla{\begin{lemm}}
\def\ela{\end{lemm}}
\def\nco{\newtheorem{coro}[axiom]{Corollary}}
\def\bco{\begin{coro}}
\def\eco{\end{coro}}
\def\ndf{\newtheorem{defn}[axiom]{Definition}}
\def\bdf{\begin{defn}}
\def\edf{\end{defn}}
\def\nrm{\newtheorem{rmrk}[axiom]{Remark}}
\def\brm{\begin{rmrk}}
\def\erm{\end{rmrk}}
\def\nex{\newtheorem{exmp}[axiom]{Example}}
\def\bex{\begin{exmp}}
\def\eex{\end{exmp}}

%% end abbreviations

%\theoremstyle{change}
%\theorembodyfont{\itshape}
\nthm
\ncn
%\theorembodyfont{\slshape}
\nla
\nco
\theoremstyle{definition}
%\theorembodyfont{\rm}
\ndf
\nrm
\nex
%% \theoremheaderfont{\scshape}

\bdoc

\title{The Existence of Cartan Connections \\ and Geometrizable Principal Bundles}
\author{Mohamed Barakat}
\thanks{I would like to thank Prof.\  E.~Ruh and Prof.\  H.~Baum for the fruitful discussions.}
\address{Lehrstuhl B f\"ur Mathematik, Technical University of Aachen, 52062 Germany}
\email{mohamed.barakat@rwth-aachen.de}

%% \date{18.2.2000}
%% Sent to Prof. Ruh:   Tue, 22 Feb 2000 14:22:05 +0100
%% Sent to Peter Quast: Mon, 15 Jan 2001 17:44:05 +0100

%\date{}

%% special abbreviations:
\def\GL{\mathrm{GL}}
\def\SO{\mathrm{SO}}
\def\g{\mathfrak{g}}
\def\h{\mathfrak{h}}
\def\p{\mathfrak{p}}
\def\Ad{\mathrm{Ad}}
\def\hor{\mathrm{hor}}
%% end abbreviations

\begin{abstract}
  The aim of this article is to give necessary and sufficient conditions
  for the existence of Cartan connections on principal bundles. For this
  purpose we recall the definition of abstract soldering forms and introduce
  the notion of geometrizable principal bundles.
\end{abstract}

\maketitle

%\rightheadtext{Cartan Connection and Geometrizable Principal Bundles}

\sec{Principal bundles and associated bundles.}
Throughout the paper let $H\to P\stackrel{\pi}{\to}M$ denote a
{\df principal $H$-bundle} over the {\df base manifold} $M$, where $P$ is
a manifold, and $H$ a Lie group with a smooth {\em free} {\em proper}
{\em right} action on $P$
\[ R:\bigg\{ \barr{rcl}
H   & \to &   \mathrm{Diff}(P); \\
h & \mapsto & R_h, \earr \]
and $M:=P/H$ the quotient manifold.

Further let $F$ be an {\df $H$-manifold}, i.e.\  a manifold  with a smooth
{\em left} $H$ action. For $p\in P$ and $\xi\in F$ denote by
\bqan p\xi:=\{(ph,h^{-1}\xi)|\ h\in H\} \eqan
the orbit of $(p,\xi)$ under the {\em natural right} action of $H$ on
$P\times F$. Then the right orbit space
\bqan E=E(P,F):=P\times_H F:=(P\times F)/H=\{p\xi|\ p\in P,\ \xi\in F\}\eqan
is called the {\df associated fiber bundle to $P$ with typical fiber $F$},
and denoted by $F\to E\stackrel{\pi_E}{\to}M$.

Every $p\in P$ with $\pi(p)=x\in M$ defines a diffeomorphism
\bqan\label{ass} p:\bigg\{ 
  \barr{rcl} F & \to & F_x:=\pi_E^{-1}(x); \\
         \xi & \mapsto & p\xi. \earr \eqan
Moreover we have: $(ph)\xi=p(h\xi)$ for every $h\in H$.

\bla\label{fund1}
  The following are in bijective correspondence:
  \ben
    \item[(i)]  A {\em global} section $f:M\to E$.
    \item[(ii)] An {\em $H$-equivariant} map $\tilde{f}:P\to F$, i.e.\
      $\tilde{f}(ph)=h^{-1}\tilde{f}(p)$ for all $h\in H$.
  \een
\ela

\bdf[Reduction]
  Let $K$ be a Lie subgroup of $H$. The principal $H$-bundle $P$ is called
  {\df $K$-reducible}, if there exists a principal $K$-bundle
  $Q\subseteq P$ over $M$, where the embedding $\iota:Q\to P$ is {\em $K$-equivariant},
  i.e.\ $\iota(qk)=\iota(q)k$ for all $q\in Q$ and $k\in K$.
  $Q$ is called a {\df $K$-reduction} of $P$.
\edf

\bdf[Symmetry breaking]
  Let $K$ be a Lie subgroup of $H$. An {\em $H$-equivariant} map
  $f:P\to H/K$ is called a {\df $K$-symmetry breaking}.
\edf

\bthm\label{fund2}
  Let $K$ be a closed Lie subgroup of $H$. The following are in pairwise
  bijective correspondence:
  \ben
    \item[(i)]   A $K$-reduction $Q\subseteq P$.
    \item[(ii)]  A {\em global} section $f:M\to P/K$.
    \item[(iii)] A $K$-symmetry breaking $\tilde{f}:P\to H/K$.
  \een
\ethm
\prf The equivalence (i) $\Leftrightarrow$ (ii) is given by factoring out the free
  action of $K$ on $Q$ and $P$. The equivalence (ii) $\Leftrightarrow$ (iii)
  is established by the fact that $P/K\cong P\times_H H/K=E(P,H/K)$ and
  Lemma \ref{fund1}.
\qed

\bdf[$H$-Structure]
  Let $M$ be an $n$-dimensional manifold and $H$ be a Lie subgroup of
  $\GL(\R^n)$. The principal bundle $H\to P\stackrel{\pi}{\to}M$ is called an
  {\df $H$-structure}, if $P$ is an $H$-reduction of the frame bundle
  $\GL(\R^n)\to L(M)\to M$.
\edf

\bdf[Pseudotensorial, Tensorial]
  Let $V$ be an $H$-module via the representation $\rho:H\to\GL(V)$. A
  $V$-valued $r$-Form $\phi:\wedge^r TP\to V$ is called
  {\df pseudotensorial} of type $(V,\rho)$ if it is $H$-equivariant, 
  i.e.\ if \[ R_h^*\phi=\rho(h^{-1})\phi,\ \ \all h\in H. \]
  A pseudotensorial $r$-form is thus an element of fixed space
  $\Gamma(\wedge^r T^*P\otimes V)^H$. \\
  $\phi$ is called {\df tensorial} if it is pseudotensorial and {\df
  horizontal}\footnote{This notion is found in \cite{kono}, Chapter II.5.
  Sharpe \cite{sha} uses `{\df semibasic}' instead.} meaning that
  \[ \phi(X_1,\ldots,X_r)=0 \] whenever at least one of the $X_i$'s is
  {\em vertical}, i.e.\  tangent to the fiber. The set of tensorial
  $r$-forms is denoted by $\Gamma_\hor(\wedge^r T^*P\otimes V)^H$.
\edf

\bs
A tensorial $r$-form factors over the induced projection
$TP\stackrel{\pi_*}{\to}TM$ and we obtain the natural identification
\bqan \Gamma_\hor(\wedge^r T^*P\otimes V)^H \cong
\Gamma(\wedge^r\pi^*T^*M\otimes V)^H.\eqan
Further if the action of $H$ on $V$ is {\em trivial}, then every tensorial
form is the pullback to the principal bundle $P$ of a $V$-valued form on
the base\footnote{This motivates the notion `{\df basic}' for such tensorial
forms, cf.\ \cite{sha}, Definition 1.5.23 and Lemma 1.5.25. See also
\cite{kono}, Example II.5.1.} manifold $M$, i.e.\
$\Gamma_\hor(\wedge^r T^*P\otimes V)^H \cong$
$\pi^*\Gamma(\wedge^r T^*M\otimes V)$.

\bdf[Bundle valued]
  Let $P$, $V$ be as above and $V\to E\stackrel{\pi_E}{\to}M$ the associated
  vector bundle. Elements of $\Gamma(\wedge^r T^*M\otimes E)$ are called
  {\df bundle valued} $r$-forms.
\edf

Now we can state the key lemma used in proof of the main theorem:
\bla\label{fund3}
  Let $P$, $V$ and $E$ be as above. The following are in bijective
  correspondence:
  \ben
    \item[(i)]  A bundle valued $r$-form on $M$.
    \item[(ii)] A tensorial $V$-valued $r$-form on $P$.
  \een
  Thus we have the natural identification:
  \bqan\label{fund3f}
    \Gamma(\wedge^r T^*M\otimes E)\cong
    \Gamma_\hor(\wedge^r T^*P\otimes V)^H. \eqan
\ela
\prf See \cite{kono}, Example II.5.2. \qed

\sec{Connections}

For the rest of the paper let $\h$ denote the Lie algebra of $H$. 

\bdf[Ehresmann connection]
  An $\h$-valued $1$-form $\gamma:TP\to\h$ is called an {\df Ehresmann
  connection} or simply a {\df connection}, if the following two conditions
  hold:
  \ben
    \item[(i)]  $R_h^*\gamma=\Ad(h^{-1})\gamma$ for all $h\in H$.
    \item[(ii)] $\gamma(X^\dagger)=X$ for all $X\in\h$.
  \een
  By $X^\dagger$ we denote the {\df fundamental vector field} of $X$ induced
  by the infinitesimal action of $\h$ on $P$ induced from the action of $H$.
\edf
Because of (i) a connection is pseudotensorial of type $(\h,\Ad)$ and
because of (ii) it is not tensorial.

\bdf[Soldering form]
  Let $M$ be an $n$-dimensional manifold, $H\to P\stackrel{\pi}{\to}M$ and
  $\rho:H\to\GL(\R^n)$ a representation. A surjective $\R^n$-valued
  tensorial $1$-form: $\theta:TP\surj\R^n$ is called a
  {\df soldering form}\footnote{In \cite{ko} $\theta$ is called a form of
  {\em soudure} (cf.\  Theorem 2 and the above lines). See also \cite{am},
  2.4, 5.1.}.
\edf

\bex[Fundamental form]
  Every $H$-structure $P$ has a natural soldering form, the so called
  {\df fundamental form}, defined by
  \[ \theta_u(X):=u^{-1}\pi_*(X)\ \ \ \all u\in P \mbox{ and } X\in TP, \]
  where $u^{-1}$ denotes the inverse of the isomorphism
  $u:\R^n\to T_{\pi(u)}M$ defined in (\ref{ass}). $\rho$ is given by the
  natural action of $H\le\GL(\R^n)$ on $\R^n$.
\eex

The definition of the fundamental form relies on the fact, that the tangent
bundle $TM$ is associated to $P$. This directly motivates the next
definition, which seems to be new.

\bdf[Geometrizable principal bundle]
  Let $M$ be an $n$-dimensional manifold. We call the principal bundle
  $H\to P\stackrel{\pi}{\to}M$ {\df geometrizable} if the tangent bundle
  $TM$ is associated to $P$, i.e.\ if there exists a representation
  $\rho:H\to \GL(\R^n)$ turning $\R^n$ into an $H$-module, such that the
  associated bundle satisfies the {\df soldering condition}:
  \bqan P\times_H\R^n\cong TM. \eqan
  We call $P$ {\df first order} geometrizable if $\rho$ is faithful,
  otherwise {\df higher order}. Once we fix a representation fulfilling the
  soldering condition, we call $P$ {\df geometrical}.
\edf

\bex
  \mbox{}
  \bit
    \item $P$ is an $H$-structure, iff $P$ is a {\df first order}
      geometrizable principal $H$-bundle. (Cf.\  Theorem \ref{main} and
      \cite{sha}, Exercise 5.3.21.)
    \item The trivial bundle $S^2\times\SO(\R^2)$ is not geometrizable.
    \item An $n$-dimensional manifold $M$ $(n\ge 3)$ is spin, iff there
      exists a geometrizable principal $\mathrm{Spin}_n(\R)$-bundle over $M$.
  \eit
\eex

It is worth mentioning, that, in contrast to relativistic gravitation
theories, most of the principal bundles appearing in classical gauge
theories are {\em not} geometrizable, unless $TM$ is trivial.

\sec{Cartan Connections}

Now we come to the main object of this paper.

\bdf[Cartan connection]
  Let $\g$ be a Lie algebra with $\h\le\g$ and $\dim\g=\dim P$. A $\g$-valued
  $1$-form $\omega:TP\to\g$ satisfying\footnote{For an alternative
  definition see \cite{am}, 2.1, 2.2.}
  \ben
    \item[(i)]   $\omega:T_p P\to\g$ is an isomorphism for all $p\in P$.
      ({\em trivialization} of $TP\to P$)
    \item[(ii)]  $R_h^*\omega=\Ad(h^{-1})\omega$ for all $h\in H$.
      ({\em pseudotensorial} of type $(\g,\Ad)$)
    \item[(iii)] $\omega(X^\dagger)=X$ for all $X\in\h$.
  \een
  We call a Cartan connection {\df reductive}, if $\h$ has an $H$-invariant
  complement\footnote{We do not require that $\p$ is a Lie subalgebra.} $\p$
  in $\g$: $\g=\h\oplus_H\p$.
\edf

In the definition of a Cartan connection we do not make use of the Lie
algebra structure of $\g$. We only use the $H$-module
structure\footnote{and the induced $\h$-module structure} of $\g$ and the
fact that $\h$ is an $H$-submodule of $\g$. The first place where we do
need the Lie algebra structure on $\g$ is in the definition of
{\em curvature} of a Cartan connection. This is exactly where the power
and flexibility of this notion lies, cf.\ \cite{ruh} and \cite{sha}.

\bthm\label{main}
  The following statements are equivalent:
  \ben
    \item[(i)] $P$ has a Cartan connection for some $\g\geq\h$.
    \item[(ii)] $P$ has a reductive Cartan connection for some $\g\geq\h$.
    \item[(iii)] $P$ has a reductive Cartan connection for some $\g\geq\h$
      with $[\p,\p]=0$.
    \item[(iv)]  $P$ has a soldering form.
    \item[(v)]   $P$ is a geometrizable bundle.
  \een
\ethm
\prf (v) $\im$ (iv): $TM$ is associated to $P$ via an $H$-module structure
  of $\R^n$. The image of the {\em identity section} $\mathrm{id}\in
  T^*M\otimes TM\cong\mathrm{End}TM$ under the bijective correspondence
  (\ref{fund3f}) for $r=1$, $E=TM$ and $V=\R^n$ (carrying the above
  $H$-module structure) is the desired soldering form. (iv) $\im$ (iii): Let
  $\theta$ be the soldering form. Due to \cite{kono}, Theorem II.2.1 there
  exists a connection form $\gamma$. $\omega=\gamma+\theta$ is a reductive
  Cartan connection with $\g:=\h\ltimes\R^n$, and $\R^n=:\p$ viewed as an
  $H$-module (and $\h$-module) via the representation $\rho$ appearing in
  the definition of the soldering form $\theta$. (iii) $\im$ (ii): Trivial.
  (ii) $\im$ (i): Trivial. (i) $\im$ (v): \cite{sha}, Theorem 5.3.15.
\qed

\ssn
Note, for a Cartan connection $\omega$ the composition
$\omega_{\g/\h}:=TP\stackrel{\omega}{\to}\g\surj\g/\h$ is a soldering form,
providing a direct proof for (i) $\im$ (iv). As mentioned in the proof,
$\omega_{\g/\h}$ is identified via (\ref{fund3f}) with the identity section
$\mathrm{id}\in\mathrm{End}TM$. Cf.\  \cite{ko}, Theorem 2, and \cite{am},
2.4, 5.1.

\bbb{aaaaaaa}
\bibitem[AM]{am} {\sc Dmitri V.~Alekseevsky and Peter W.~Michor.} {\sl
Differential geometry of Cartan connections} Publ.\  Math.\  Debrecen
{\bf 47} (1995), 349-375. Preprint
\href{http://xxx.uni-augsburg.de/abs/math/9412232}{math.DG/9412232}.

\bibitem[Ko]{ko} {\sc S.~Kobayashi.} {\sl On Connections of Cartan.}
Canad.\  J.\  Math.\  {\bf 8} (1956), 145-156.

\bibitem[KoNo]{kono} {\sc S.~Kobayashi and K.~Nomizu.} {\sl Foundations of
differential geometry.} Volume I.\  Interscience Publishers, New York, 1963.

\bibitem[Ruh]{ruh} {\sc E.~Ruh.} {\sl Almost flat manifolds.} J.\  of
Differential Geom.\  {\bf 17} (1982), 1-14.

\bibitem[Sh]{sha} {\sc R.W.~Sharpe.} {\sl Differential Geometry.}
Cartan's Generalization of Klein's Erlangen Program. Springer 1996.

\ebb

\edoc